\documentclass[12pt]{amsart}
\usepackage{amsthm,amssymb}
\newcommand{\ol}{\overline}

\newcommand{\iy}{\infty}

\newcommand{\sg}{\sigma}

\newcommand{\bT}{\mathbb{T}}

\newtheorem{theorem}{\bf  Theorem}

\theoremstyle{remark}

\renewenvironment{proof}{ P\,r\,o\,o\,f.}{$\Box$}

\begin{document}
\begin{center}
    {\bf A SIMPLE PROOF OF THE MATRIX-VALUED FEJ\'{E}R-RIESZ THEOREM}

\vskip+0.5cm
   {Lasha Ephremidze, \fbox{Gigla Janashia} and Edem Lagvilava  }
\end{center}

\footnotetext { Received by the editors August 15, 2007.



}

\vskip+0.5cm

{\small {\bf Abstract.} A very short proof of the Fej\'{e}r-Riesz
lemma is presented in the matrix case.}

\vskip+0.5cm

The following fundamental result in matrix spectral factorization
theory belongs to Wiener [8] (see also [5], [4]):
\smallskip

{\bf Theorem.} {\em Let
\begin{equation}
    S(z)\sim \sum_{n=-\iy}^\iy\sg_nz^n,
\end{equation}
$|z|=1$,  $\sg_k$ are $r\!\times\!r$ matrix coefficients, be a
positive definite matrix-function with integrable entries, $S(z)\in
L_1(\bT)$. If the logarithm of the determinant is integrable, $\log
\det S(z)\in L_1({\mathbb T})$, then there exists a factorization
\begin{equation}
    S(z)=\chi^+(z)(\chi^+(z))^*,
\end{equation}
where
\begin{equation}
    \chi^+(z)= \sum_{n=0}^\iy\rho_nz^n,
\end{equation}
$|z|<1$, is an analytic matrix-function with entries from the Hardy
space $H_2$, $\chi^+(z)\in H_2$, and the determinant of which is an
outer function.}
\smallskip

The relation (2) is assumed to hold a.e. on the unit circle $\bT$
and $(\chi^+)^*=(\ol{\chi^+})^T$ is the adjoint of $\chi^+$.

It is well-known that if $S(z)$ in (1) is a Laurent polynomial of
order $m$, then $\chi^+(z)$ in (3) is a polynomial of the same order
$m$. This result is known as the Fej\'{e}r-Riesz lemma in the scalar
case and it was  generalized by M. Rosenblatt [7] and Helson [4] to
the matrix case using the linear prediction theory of
multidimensional weakly stationary processes as in the proof of the
existence theorem above. One can find a constructive but long proof
of the matrix Fej\'{e}r-Riesz lemma in [3] as well. Below, we
present a simple, transparent and natural proof of this assertion.

\begin{theorem}
If $S(z)=\sum_{n=-m}^m\sg_nz^n$, then $\chi^+(z)=
\sum_{n=0}^m\rho_nz^n$.
\end{theorem}

We use the generalization of Smirnov's theorem (see [6], p. 109)
which claims that if the boundary values of an analytic function
$f(z)=g(z)/h(z)$, where $g\in H_{p_1}$ and $h$ is an outer function
from $H_{p_2}$, belongs to $L_p({\mathbb T})$, then $f\in H_p$\,.

\smallskip

\begin{proof}
We know that
\begin{equation}
(\chi^+(z))^*\in L_2^-(\bT)
\end{equation}
in general, and it suffices to show that
\begin{equation}
    z^m(\chi^+(z))^*\in L_2^+(\bT)=H_2,
\end{equation}
where $L_2^-(\bT)$ and $L_2^+(\bT)$ are the classes of square
integrable functions with, respectively, positive and negative
Fourier coefficients equal to $0$, and the latter is naturally
identified with $H_2$.

It follows from (2) that
\begin{equation}
 (\chi^+(z))^{-1}z^m S(z)=z^m(\chi^+(z))^*
\end{equation}
for a.a. $z\in\bT$. The matrix-function $$(\chi^+(z))^{-1}=\frac
1{\det\chi^+(z)}A(z)$$ is analytic in the unit circle, where
$A(z)\in H_{2/r}$\,. Consequently, since  $z^mS(z)\in H_\iy\,$ by
hypothesis, the entries of the left-hand side matrix in (6) can be
represented as the ratios of some functions from $H_{2/r}$ and the
outer function $\det\chi^+(z)\in H_{2/r}$, while their boundary
values belong to $L_2(\bT)$ because of (6) and (4). Thus, by virtue
of the above mentioned generalization of Smirnov's theorem,
$(\chi^+(z))^{-1}z^m S(z)\in H_2=L_2^+(\bT)$ and (5) follows again
from (6).
\end{proof}
\smallskip

The same  idea can be used to prove the uniqueness (up to a constant
unitary matrix) of the spectral factorization (2). Indeed, assume
$S(z)=\chi^+_1(z)(\chi^+_1(z))^*$ together with (2), where
$\chi^+_1(z)\in H_2$ and $\det\chi^+_1(z)$ is outer. Then
\begin{equation}
 (\chi^+(z))^{-1}\chi^+_1(z)((\chi^+(z))^{-1}\chi^+_1(z))^*=I,
\end{equation}
so that the analytic matrix-function
$U(z):=(\chi^+(z))^{-1}\chi^+_1(z)$, $|z|<1$, is unitary  for a.a.
$z\in\bT$. Thus, the boundary values of $U(z)$ belongs to $L_\iy$
and, as in the proof of Theorem 1, we have $U(z)\in H_\iy$. By
changing the roles of $\chi^+$ and $\chi^+_1$ in this discussion, we
get $(\chi^+_1(z))^{-1}\chi^+(z)\in H_\iy$. But
$(U(z))^*=(\chi^+_1(z))^{-1}\chi^+(z)$ for a.a. $z\in\bT$, by virtue
of (7). Thus the boundary values of $U(z)$ as well as its conjugate
belongs to $L_\iy^+(\bT)$ which implies that $U(z)$ is constant.

A further development of the circle of ideas presented in this paper
leads to the constructive analytic proof of the existence theorem,
formulated in the beginning, as well as to the efficient algorithm
for approximate computation of the spectral factor (3) for a given
matrix-function (1) (see [2], [1]). Such an algorithm is very
important for practical applications and it has been searched since
Wiener's period.

\vskip+0.5cm

{\small {\sc A. Razmadze Mathematical Institute,  Tbilisi 0193,
Georgia}

{\em  E-mail addresses}: lasha@rmi.acnet.ge; edem@rmi.acnet.ge

\vskip+0.2cm

Current address of L. Ephremidze: {\sc Department of Mathematics,
Tokai University, Shizuoka,   424-8610, Japan}

 {\em E-mail address}: 05k054@scc.u-tokai.ac.jp
}

\end{document}